\documentclass[a4paper]{article}
\usepackage[english]{babel}  
\usepackage[T1]{fontenc}  
\usepackage{amsmath,amssymb}  
\usepackage{color}
\usepackage{geometry}
\usepackage[all]{xy}
\usepackage{amsthm}
\usepackage{varioref}
\usepackage{hyperref}
\usepackage{graphicx}

\hypersetup{
						colorlinks=true,
						breaklinks=true,
						urlcolor= blue,
						linkcolor= black,
						anchorcolor = black,
						citecolor = black,
						filecolor = black,
					}
						
\newcommand{\C}{ {\mathbb C} }
\newcommand{\R}{ {\mathbb R} }
\newcommand{\D}{ {\mathbb D} }
\newcommand{\Z}{ {\mathbb Z} }
\newcommand{\N}{ {\mathbb N} }
\newcommand{\Q}{ {\mathbb Q} }
\newcommand{\PP}{ {\mathbb P} }

\newcommand{\cP}{{\mathcal P}}
\newcommand{\cE}{{\mathcal E}}
\newcommand{\cA}{{\mathcal A}}
\newcommand{\cZ}{{\mathcal Z}}
\newcommand{\cV}{{\mathcal V}}
\newcommand{\cS}{{\mathcal S}}
\newcommand{\cN}{{\mathcal N}}

\newtheoremstyle{thm}
  {0pt}
  {5pt}
  {}
  {}
  {\bfseries}
  {:}
  {.5em}
  {}
  
\theoremstyle{thm}
\newtheorem{thm}{Theorem}
\labelformat{thm}{theorem~#1}
\newtheorem{cor}{Corollary}
\labelformat{cor}{corollary~#1}
\newtheorem{prop}{Proposition}[section]
\labelformat{prop}{proposition~#1}

\newtheorem*{propnonnum}{Proposition}
\newtheorem*{defn}{Definition}
\newtheorem*{demo}{Proof}
\newtheorem{exempl}{Example}
\labelformat{exempl}{example~#1}

\newtheoremstyle{remarque}
  {3pt}
  {5pt}
  {}
  {}
  {\itshape}
  {:}
  {.5em}
  {}

\theoremstyle{remarque}
\newtheorem*{rmq}{Remark}

\date{}

\begin{document}

\title{LVMB manifolds and simplicial spheres}
\author{Jérôme Tambour}

\maketitle

\begin{abstract}
\noindent
LVM and LVMB manifolds are a large family of examples of non kähler manifolds. For instance, Hopf manifolds and Calabi-Eckmann manifolds can be seen as LVMB manifolds. The LVM manifolds have a very natural action of a real torus and the quotient of this action is a simple polytope. This quotient allows us to relate closely LVM manifolds to the moment-angle manifolds studied (for example) by Buchstaber and Panov. Our aim is to generalize the polytopes associated to LVM manifolds to the LVMB case and study the properties of this generalization. In particular, we show that the object obtained belongs to a very large class of simplicial spheres. Moreover, we show that for every sphere belonging to this class, we can construct a LVMB manifold whose associated sphere is the given sphere. We use this latter result to show that many moment-angle complexes can be endowed with a complex structure (up to product with circles). 
\end{abstract}

\maketitle

\section*{Introduction}

It is not easy to construct non kähler compact complex manifolds. The simplest example is the well known Hopf manifold (\cite{H}, 1948), which gives a complex structure on the product of spheres $S^{2n+1}\times S^1$ as a quotient of $\C^n\backslash\{0\}$ by the action of a discrete group. The Hopf manifold has many generalizations: firstly, by Calabi and Eckmann \cite{CE} who give a structure of complex manifold on any product of spheres (of odd dimension). Then by Santiago Lopez de Medrano, Alberto Verjovsky (\cite{LdM} and \cite{LdMV}) and Laurent Meersseman \cite{M}. In these last generalizations, the authors obtain complex structures on products of spheres, and on connected sums of products of spheres, also constructed as a quotient of an open subset in $\C^n$ but by the action of a non discrete group. These manifolds are known as \emph{LVM manifolds}.

\vskip 4mm

The construction in \cite{M} has (at least) two interesting features: on the one hand, the LVM manifolds are endowed with an action of the torus $(S^1)^n$ whose quotient is a simple convex polytope and the combinatorial type of this polytope characterizes the topology of the manifold. On the other hand, for every simple polytope $P$, it is possible to construct a LVM manifold whose quotient is $P$.

\vskip 4mm

In \cite{B}, Frédéric Bosio generalizes the construction of \cite{M} emphasing on the combinatorial aspects of LVM manifolds. This aim of this paper is to study the \emph{LVMB manifolds} (i.e. manifolds constructed as in \cite{B}) from the topological and combinatorial viewpoints. In particular, we will generalize the associated polytope of a LVM manifold to our case and prove that this generalization belongs to a large class of simplicial spheres (named here \emph{rationally starshaped spheres}). 

\vskip 4mm

In the first part, we briefly recall the construction of the LVMB manifolds as a quotient of an open set in $\C^n$ by a holomorphic action of $\C^m$. In the second part, we study \emph{fundamental sets}, the combinatorial data describing a LVMB manifold and their connection to pure simplicial complexes. Mainly, we show that an important property appearing in \cite{B} (the $SEU$ property) is related to a well-known class of simplicial complexes: the pseudo-manifolds. In the same part, we also introduce the simplicial complex associated to a LVMB manifold and we show that this complex generalizes the associated polytope of a LVM manifold. In the third and forth parts, we are mainly interested in the properties of this complex and we show that the complex is indeed a simplicial sphere. To do that, we have to study another action whose quotient is a toric variety closely related to our LVMB manifold. This action was already studied in \cite{MV} and \cite{CFZ} but we need a more thorough study. Finally, in the fifth part, we make the inverse construction: starting with a rationally starshaped sphere, we construct a LVMB manifold whose associated complex is the given sphere. Using this construction, we show an important property for moment-angle complexes: up to a product of circles, every moment-angle complex arising from a starshaped sphere can be endowed with a complex structure of LVMB manifold.

\vskip 4mm 

To sum up, we prove the following theorems:

\begin{thm} \label{thmprinc} Let $\cN$ be a LVMB manifold. Then its associated complex $\cP$ is a rationally starshaped sphere. Moreover, if $\cN$ is a LVM manifold, then $\cP$ can be identified with (the dual of) its associated polytope.
\end{thm}

\vskip 2mm

\begin{propnonnum}
Every rationally starshaped sphere can be realized as the associated complex $\cP$ for some LVMB manifold.
\end{propnonnum}

\vskip 2mm

\begin{thm} Up to a product of circles, every moment-angle complex arising from a starshaped sphere can be endowed with a complex structure of LVMB manifold.
\end{thm}


\subsection*{Notations}
In this short section, we fix several notations which will be used throughout the text:

\begin{itemize}
	\item Si $A$ is a subset of a set $V$, we denote $V\backslash A$ its complement in $V$, or simply $A^c$ if no confusion can be made.
	\item We put $I_z=\{\ k\in\{1,\dots ,n\}\ /\ z_k\neq 0\ \}$ for every $z$ in $\C^n$.
	\item $\D$ is the closed unit disk in $\C$ and $S^1$ its boundary.
	\item Moreover, $exp$ will be the map $\C^n\rightarrow (\C^*)^n$ defined by 
		\[
		exp(z)=\left(e^{z_1},\dots ,e^{z_n}\right)
		\] 
		(where $e$ is the usual exponential map of $\C$).
	\item If $m\in\Z^n$, we will denote $X_n^m$ the character of $\left(\C^*\right)^n$ defined by $ X_n^m(z)=z_1^{m_1}\dots z_n^{m_n}$. 
	\item And the one-parameter subgroups of $(\C^*)^n$ will be denoted $\lambda_n^m$: 
		\[
		\lambda_n^m(t)=(t^{m_1},\dots,t^{m_n})
		\]
	\item $<,>$ is the usual \textbf{non hermitian} inner product on $\C^n:$ 
		\[
		<z,w>=\sum_{j=1}^nz_jw_j
		\]
	\item We will identify $\C^m$, as a $\R$-vector space, to $\R^{2m}$ via the morphism 
		\[
		z\mapsto \left(Re\left(z_1\right),\dots,Re\left(z_n\right),Im\left(z_1\right),\dots,Im\left(z_n\right)\right)
		\]
	\item We set $Re\left(z\right)=\left(Re\left(z_1\right),\dots,Re\left(z_n\right)\right)$
	\item As well, we set $Im\left(z\right)=\left(Im\left(z_1\right),\dots,Im\left(z_n\right)\right)$, so 
		\[
		z=Re\left(z\right)+iIm\left(z\right)=\left(Re\left(z\right),Im\left(z\right)\right)
		\]
	\item In $\R^n$, Conv(A) is the convex hull of a subset A. 
	\item Si $A$ is a nonempty subset of $\R^n$, the set of all nonnegative linear combinations 
		\[
		x=\displaystyle{\sum_{j=1}^k\lambda_ja_j},\hskip 2mm k\in\N^*,\hskip 2mm \lambda_1\geq0,\dots,\lambda_k\geq0,\hskip 2mm a_1,\dots,a_k\in A 
		\] 
\noindent		 
of elements of $A$ is called the \emph{positive hull} of $A$ and denoted $pos(A)$. If $A=\emptyset$, we define $pos(\emptyset)=\{0\}$. 
	\item For every $v$ in $\R^n$, we set $\widetilde{v}=(1,v)\in\R^{n+1}$.
\end{itemize}

\vskip 4mm 
\section{Construction of the LVMB manifolds}

In this section, we briefly recall the construction of LVMB manifolds, following the notation of \cite{B}. Let $M$ and $n$ be two positive integers such that $n\geq M$. A \emph{fundamental set} is a nonempty set $\cE$ consisting of subsets of $\{1,\dots,n\}$ having $M$ elements. Elements of $\cE$ are called \emph{fundamental subsets}. 

\vskip 2mm

\begin{rmq} For practical reasons, we sometimes consider fundamental sets whose elements do not belong to $\{1,\dots,n\}$ but to another finite set with $n$ elements (usually $\{0,1,\dots, n-1\}$).
\end{rmq}

\vskip 4mm

Let $P$ be a subset of $\{1,\dots,n\}$. We say that $P$ is \emph{acceptable} if $P$ contains a fundamental subset. We define $\cA$ as the set of all acceptable subsets. Finally, an element of $\{1,\dots,n\}$ will be called \emph{indispensable} if it belongs to every fundamental subset of $\cE$. We say that $\cE$ is of type $(M,n)$ (or $(M,n,k)$ if we want to emphasize the number $k$ of indispensable elements).

\vskip 4mm

\begin{exempl}\label{exemple} For instance, \[\cE=\{\ \{1,2,5\},\{1,4,5\},\{2,3,5\},\{3,4,5\}\ \}\] is a fundamental set of type $(3,5,1)$.
\end{exempl}

\vskip 4mm

Two combinatorial properties (named $SE$ and $SEU$\footnote{for Substitute's Existence and Substitute's Existence and Uniqueness} respectively, cf. \cite{B}) will be very important in the sequel:

\vskip 2mm

\[
(SE) \hskip 5mm \forall\ P\in\cE,\hskip 2mm \forall\ k \in\{1,\dots,n\},\hskip 2mm \exists\ k'\in P;\hskip 2mm (P\backslash\{k'\})\cup\{k\}\in\cE
\]

\vskip 4mm

\[
(SEU) \hskip 5mm \forall\ P\in \cE, \hskip 2mm \forall\ k\in\{1,\dots,n\}, \hskip 2mm \exists!\ k'\in P; \hskip 2mm (P\backslash\{k'\})\cup\{k\}\in\cE
\]

\vskip 7mm

Moreover, a fundamental set $\cE$ is \emph{minimal} for the $SEU$ property if it verifies this property and has no proper subset $\widetilde\cE$ such that $\widetilde\cE$ verifies the $SEU$ property. Finally, we associate to $\cE$ two open subset $\cS$ in $\C^n$ and $\cV$ in $\PP^{n-1}$ (where $\PP^{n-1}$ is the complex projective space of dimension $n-1$) defined as follows: 

\vskip 2mm

\label{ouvertS}
\[
\cS=\{\ z\in \C^n/\  I_z\in \cA\ \}
\]

\vskip 4mm

The open subset $\cS$ is the complement of an arrangement of coordinate subspaces in $\C^n$. Indeed, we have the following description of $\cS$:

\vskip 2mm

\begin{prop} We have 

\[
\displaystyle{\cS=\C^n\ \backslash\ \bigcup_{(i_1,\dots,i_k)\notin\cP} \left\{\ z\ \middle/\ z_{i_1}=\dots=z_{i_k}=0 \right\}\ }
\]

\end{prop}

\vskip 4mm

\begin{demo} Let $z$ be an element of $\C^n$. If $\sigma=(i_1,\dots,i_k)\subset\{1,\dots,n\}$, we denotes $L_\sigma$ for the coordinate subspace 

\[
\displaystyle{L_\sigma= \left\{\ z\ \middle/\ z_{i_1}=\dots=z_{i_k}=0 \right\}\ }
\]

\vskip 4mm

The elements of $L_\sigma$ are exactly the elements $z$ such that $I_z^c$ contains $\sigma$. If $z$ belongs to $\cS$, then, by definition, $I_z^c$ is a face of $\cP$. So, if $z$ belongs to some $L_\sigma$, we have $\sigma\in\cP$. Conversely, if $z$ is not an element of $L_\sigma$ for every $\sigma\notin\cP$, then, since $z$ belongs $L_{I_z^c}$, we conclude that $I_z^c$ is a face of $\cP$, that is, $z$ belongs to $\cS$.
\end{demo}

\vskip 4mm

We also remark that, for all $t\in\C^*$, we have $I_{tz}=I_z$ so the following definition is consistent:

\vskip 2mm

\[
\cV=\{\ [z]\in \PP^{n-1}/\  I_z\in \cA\ \}
\]

\vskip 5mm

\begin{exempl} In Example \ref{exemple}, $\cE$ verifies the $SEU$ property and the set $\cS$ is 

\[
\cS=\left\{\ \left(z_1,z_2,z_3,z_4,z_5\right)\ \middle/\ \left(z_1,z_3\right)\neq0,\ \left(z_2,z_4\right)\neq0,\ z_5\neq0\ \right\}
\]

\vskip 2mm

so 

\[
\cS \simeq\ \left(\C^2\backslash\{0\}\right)^2\times\C^* 
\]

\end{exempl}

\vskip 5mm

Now, we suppose that $M=2m+1$ is odd and we fix $l=(l_1,\dots,l_n)$ a family of elements of $\C^m$. We can define an action (called \emph{acceptable holomorphic action}) of $\C^*\times\C^m$ on $\C^n$ defined by:

\[
(\alpha,T)\cdot z=\left(\ \alpha e^{<l_1,T>}z_1,\ \dots,\ \alpha e^{<l_n,T>}z_n\ \right)\hskip 4mm \forall\ (\alpha,T,z) \in \C^*\times\C^m\times\C^n
\]

\vskip 4mm

\begin{rmq} If $m=0$, the previous action is just the classical one defining $\PP^{n-1}$. 
\end{rmq}

\vskip 4mm

In the sequel, we focus only on families $l$ such that for every $P\in \cE$, $(l_p)_{p\in P}$ spans $\C^m$ (seen as a $\R$-affine space). In this case, we say that $(\cE,l)$ is a \emph{acceptable system}. If $\alpha\in\C^*$, $T\in\C^m$ and $z\in\C^n$, we have $I_{(\alpha,T).z}=I_z$, so $\cS$ is invariant for the acceptable holomorphic action. Moreover, the restriction of this action to $\cS$ is free (cf. \cite{B}, p.1264).

\vskip 4mm

As a consequence, we denote $\cN$ the orbit space for the action of $\C^*\times\C^m$ restricted to $\cS$. Notice that we can also consider an action of $\C^m$ on $\cV$ whose quotient is again $\cN$. We call $(\cE,l)$ a \emph{good system} if $\cN$ is compact and can be endowed with a complex structure such that the natural projection $ \cS\rightarrow{\cN}$ is holomorphic. Such a manifold is known as a \emph{LVMB manifold}. Since a quotient of a holomorphic manifold by a free and proper action (cf. \cite{Hu}, p.60) can be endowed with a complex structure, we only have to check whenever the action is proper. Here, according to \cite{Hu},\ p.59, we define a proper action of a Lie Group $G$ on a topological space $X$ as a continuous action such that the map $G \times X \rightarrow X \times X$ defined by $(g,x)\ \mapsto (g\cdot x,x)$ is proper. Notice that in our case, the group $G$ is not discrete.  

\vskip 4mm

We recall the following definition from \cite{B}:   

\vskip 2mm

\begin{defn} Let $(\cE,l)$ be a acceptable system. We say that it verifies the \emph{imbrication condition} if for every $P,Q$ in $\cE$, the interiors of the convex hulls $Conv(l_p,p\in P)$ and $Conv(l_q,q\in Q)$ have a common point.
\end{defn}

\vskip 4mm

In \cite{B}, the following fundamental theorem is proved:

\vskip 2mm

\begin{thm}[\cite{B}, p.1268] \label{thmBosio1} An acceptable system is good if and only if $(\cE,l)$ verifies the $SE$ property and the imbrication condition.
\end{thm}

\vskip 2mm

\begin{rmq} In \cite{B}, it is also proved that a good system is minimal for the $SEU$ property.
\end{rmq}

\vskip 4mm

Finally, a \emph{LVM manifold} is a manifold constructed as in \cite{LdMV} or \cite{M}. We don't explain here the whole construction of the LVM manifolds. The only thing we need here is that is a special case of LVMB. Indeed, we have the following theorem (and we use this theorem as a definition of LVM manifolds):

\vskip 2mm

\begin{thm}[\cite{B}, p.1265] \label{LVM} Every LVM manifold is a LVMB manifold. To be more precise, let $\mathcal O$ be the set of points of $\C^m$ which are not in the convex hull of any subset of $l$ with cardinal $2m$. Then, a good system $(\cE,l)$ is the good system of a LVM manifold if and only if there exists an bounded component $O$ in $\mathcal O$ such that $\cE$ is exactly the set of subsets $P$ of $\{1,\dots,n\}$ with $(2m+1)$ for cardinal such that $O$ is contained in $Conv(l_p,p\in P)$.
\end{thm}

\vskip 2mm

\begin{rmq} In particular, a LVM manifold is a compact complex manifold.
\end{rmq}

\vskip 4mm

\begin{exempl} We come back to example \ref{exemple}. If we set $l_1=l_3=1$, $l_2=l_4=i$ and $l_5=0$, then the imbrication condition is fulfilled and $(\cE,l)$ is a good system. As a consequence, Theorem \ref{thmBosio1} and Theorem \ref{LVM} imply that $\cN$ can be endowed with a structure of a LVM manifold. In \cite{LdMV}, the LVM manifolds constructed from a good system of type $(3,n,k)$ are classified up to diffeomorphism. Here, $\cE$ has type $(3,5,1)$ and we have $\cN\approx S^3\times S^3$. Notice that we can also use the theory of moment-angle complex (cf. the section \ref{cma}) to do this calculation.  
\end{exempl}

\vskip 4mm

To conclude this section, we recall how to construct the polytope associated to a LVM manifold $(\cE,l)$. It is clear that the natural action of $(S^1)^n$ on $\C^n$ preserves $\cS$ and that this action commutes with the given holomorphic action. So, we have an induced action of $(S^1)^n$ on $\cN$. Up to a translation of the $l_j$ (which does not change the quotient $\cN$), Theorem \ref{LVM} allow us to assume that $0$ belongs to $Conv(l_1,\dots,l_n)$ (this condition is known as \emph{Siegel's condition}) and in this case, the quotient ${\mathfrak P}$ of the action of $(S^1)^n$ on $\cN$ can be identified with:

\[
{\mathfrak P}=\left\{\ (r_1,\dots,r_n)\in(\R_+)^n\ \middle/\ \sum_{i=1}^n r_jl_j=0,\ \sum_{i=1}^n r_j=1\ \right\}
\]  

\vskip 4mm

Since a LVM manifold is compact, this set is clearly a polytope and it can be shown that it is simple (cf. \cite{BM}, Lemma $0.12$). The polytope $P$ is called the \emph{polytope associated} to the LVM manifold $\cN$.

\vskip 4mm

\begin{exempl} For the previous example, we have to perform a translation on the vectors $l_j$ with the view to respect the Siegel condition. For example, $\cN$ is also the LVM manifold associated with the good system $(\cE,\lambda)$ where $\lambda_1=\lambda_3=\frac{3}{4}-\frac{i}{4}$, $\lambda_2=\lambda_4=-\frac{1}{4}+\frac{3i}{4}$ and $\lambda_5=-\frac{1}{4}-\frac{i}{4}$. A calculation shows that the polytope ${\mathfrak P}$ is the square 

\[
{\mathfrak P}=\displaystyle{\left\{\ \left(\frac{1}{4}-r_3,\frac{1}{4}-r_4,r_3,r_4,\frac{1}{2}\right)\ \middle/\ r_3,r_4\in \left[0,\frac{1}{4}\right]\ \right\}}
\] 

\end{exempl}

\section{Fundamental sets and associated complex}

In this section, we will briefly study the fundamental sets introduced in \cite{B} and recalled above, and construct a simplicial complex whose combinatorial properties reflect the geometry of a LVMB manifold. Here, $\cE$ is a fundamental set of type ($M,n$). The integer $M$ is not supposed to be odd. For the moment, our aim is to relate the above properties to some classical ones of simplicial complexes. 

\vskip 4mm

Let us begin with some terminology. \emph{Faces}, or \emph{simplices} are subsets of a simplicial complex. If a complex $K$ is \emph{pure-dimensional} (or simply \emph{pure}), the simplices of maximal dimension $d$ are named \emph{facets} and the faces of dimension 0 (resp. $d-1$) are the \emph{vertices} (resp. the \emph{ridges}) of the complex.

\vskip 4mm

The first important definition in this paper is the following:

\vskip 2mm

\begin{defn} Let $\cE$ be a fundamental set of type $(M,n)$. The \emph{associated complex} of $\cE$ is the set $\cP$ of subsets of $\{1,\dots,n\}$ whose complement (in $\{1,\dots,n\}$) is acceptable.
\end{defn}

As we will see later, this complex is the best choice for a combinatorial generalization of the associated polytope of a LVM manifold.

\vskip 4mm

First properties of $\cP$ are the following:

\vskip 2mm

\begin{prop} \label{pcomplex} Let $\cE$ be a fundamental set of type $(M,n,k)$. Then, its associated complex $\cP$ is a simplicial complex on $\{1,\dots,n\}$. Moreover, $\cP$ is pure-dimensional of dimension $(n-M-1)$ and has $(n-k)$ vertices. These vertices are precisely the non-indispensable elements of $\{1,\dots,n\}$ for $\cE$ and the facets are exactly the complements of the subsets of $\cE$.
\end{prop}

\begin{demo} $\cP$ is obviously a simplicial complex. Moreover, the maximal simplices of $\cP$ are exactly the complements of minimal subsets of $\cA$, i.e the fundamental subsets. The latters have the same number $M$ of elements, so every maximal simplex of $\cP$ has $n-M$ elements. Finally, an element $j\in\{1,\dots,n\}$ is a vertex of $\cP$ if and only if $\{1,\dots,n\}\backslash\{j\}$ contains a fundamental subset, that is, $j$ is not indispensable.  
\end{demo}

\vskip 4mm

\begin{exempl} The complex $\cP$ associated to the fundamental set of Example \ref{exemple} is the complex with facets $\{ \{1,2\},\{2,3\},\{3,4\},\{1,4\}\}$. So, $\cP$ is the boundary of a square.
\end{exempl}

\vskip 2mm

\begin{rmq} Conversely, every pure complex can be realized as the associated complex of a fundamental set: let $\cP$ be a pure-dimensional simplicial complex on the set $\{1,\dots,v\}$ with dimension $d$ and $v$ vertices. Then, for every integer $k$, there exists  two integers $M,n$ and a fundamental set $\cE$ of type $(M,n,k)$ whose associated complex is $\cP$. If $k$ is fixed, this fundamental set is unique.
\end{rmq}

\vskip 4mm

Moreover, the $SEU$ property can be expressed as a combinatorial property of $\cP$:

\vskip 2mm

\begin{prop} \label{facettes} $\cE$ satisfies the $SEU$ property if and only if 

\[
\forall\ Q\in \cP_{max},\ \forall\ k\in \{1,\dots,n\},\ \exists!\ k'\notin Q;\ (Q\cup\{k'\})\backslash\{k\}\in\cP_{max} 
\]

where $\cP_{max}$ is the set of facets of $\cP$.
\end{prop}

\vskip 2mm

\begin{demo} First, we assume that $\cE$ verifies the $SEU$ property. Let $Q$ be a facet of $\cP$ and $k\in\{1,\dots,n\}$. Then $P=Q^c$ belongs to $\cE$. The $SEU$ property implies that there is $k'$ in $P$ (so $k'\notin Q$) such that $P'=(P\backslash\{k'\})\cup\{k\}$ belongs to $\cE$. As  claimed, $P'^c$ is also a facet of $\cP$. Moreover, we obviously have $P'^c=(Q\cup\{k'\})\backslash\{k\}$. Finally, if  $Q''=(Q\cup\{k''\})\backslash\{k\}$ is a facet of $\cP$ with $k''\neq k'$ and $k''\notin Q$, then $(P\backslash\{k''\})\cup\{k\}$ is in $\cE$, which contradicts the $SEU$ property. The proof of the converse is analogous.  
\end{demo}

\vskip 4mm

\begin{cor} \label{crete} Let $\cE$ be a fundamental set. Then its associated complex $\cP$ satisfies the $SEU$ property if and only if every ridge of $\cP$ is contained in exactly two facets of $\cP$.
\end{cor}

\begin{demo} To begin, we assume that $\cE$ verifies the $SEU$ property. Let $Q$ be a ridge of $\cP$. By definition, $Q$ is included in a facet $P$ of $\cP$. We put $P=Q\cup\{k\}, k\in\{1,\dots,n\}\backslash Q$. By proposition \ref{facettes}, there exists $k'\notin P$ (and so $k\neq k'$) such that $P'=(P\cup\{k'\})\backslash\{k\}$ is a facet of $\cP$. We have $P'=Q\cup\{k'\}$ so $Q$ is contained in at least two facets of $\cP$. Let assume that $Q$ is contained in a third facet $P''=Q\cup\{k''\}$. In this case, we have $P''=(P\cup\{k''\})\backslash\{k\}$, which contradicts proposition \ref{facettes}.
\\Conversely, let $Q$ be a facet of $\cP$ and $k\in\{1,\cdots,n\}$. If $k\in Q$, then $P=Q\backslash\{k\}$ is a ridge of $\cP$ and by hypothesis, $P$ is  contained in exactly two facets $Q_1$ and $Q_2$. One of them, say $Q_1$, is $Q$. The other is $Q_2$ and we have $Q_2=P\cup\{k'\}$. Then we have $k'\notin Q$ (on the contrary, we would have $Q_2=Q=Q_1$) and $Q_2=(Q\cup\{k'\})\backslash\{k\}$. Moreover, if $Q_3=(Q\cup\{k''\})\backslash\{k\}$ is another facet of $\cP$ with $k''\notin Q$, then $Q_3$ contains $P$ and by hypothesis, $Q_3=Q_2$ (i.e $k''=k'$). If $k\notin Q$, we remark that the element $k'\in\{1,\dots,n\}$ such that $Q'=(Q\cup\{k'\})\backslash\{k\}$ is a facet of $\cP$ is $k'=k$. Indeed, if $k'=k$, then $Q'=Q$ is a facet of $\cP$. And if $k'\neq k$, then $k\notin Q\cup\{k'\}$ and, as a consequence, $Q'=Q\cup\{k'\}$ is not in $\cP$.  
\end{demo}

\vskip 4mm

\begin{defn} Let $\cE$ be a fundamental set of type $(M,n)$. We define the (unoriented) graph $\Gamma$ by requiring that its vertices are fundamental subsets of $\cE$ and two vertices $P$ and $Q$ are related by an edge  if and only if there exist $k\notin P,k'\in P$ such that $Q=(P\backslash\{k'\})\cup\{k\}$. Equivalently, we relate two subsets of $\cE$ if and only if they differ exactly by one element. $\Gamma$ is called the \emph{replacement graph} of $\cE$.
\end{defn}

\vskip 4mm

\begin{prop} Let $\cE$ be a fundamental set of type $(M,n)$ which verifies the $SEU$ property. Then, there exist an integer $p\in\N^*$, and fundamental sets $\cE_j$ of type $(M,n)$ which are minimal for the $SEU$ property and such that $\cE$ is the disjoint union $\displaystyle{\bigsqcup_{j=1}^p\cE_j}$.
\end{prop}

\begin{demo} We proceed by induction on the cardinal of $\cE$. If $\cE$ is minimal for the $SEU$ property, then there is nothing to do. Let assume that it is not the case: there exists a proper subset $\cE_1$ of $\cE$ which is minimal for the $SEU$ property. We put $\overline\cE$ for its complement $\cE\backslash\cE_1$. It is obvious that $\overline\cE$ is a fundamental set (of type $(M,n)$). We claim that $\overline\cE$ verifies the $SEU$ property. Let $P$ be an element of $\overline\cE$ and $k\in\{1,\dots,n\}$. If $k\in P$, then, putting $k'=k$, we have that $(P\backslash\{k'\})\cup\{k\}=P$ is an element of $\overline\cE$. It is the only choice (for $k'$) since $P$ is in $\cE$ and $\cE$ verifies the $SEU$ property. Let assume now that $k$ is not in $P$. Since $P$ is an element of $\cE$, there exists exactly one $k'\in P$ such that $P'=(P\backslash \{k'\})\cup\{k\}$ is an element of $\cE$, too. We claim that $P'$ cannot be in $\cE_1$. Indeed, if it were the case, since $\cE_1$ is minimal for the $SEU$ property, there would exist exactly one $k''\in P'$ such that $P''=(P'\backslash\{k''\})\cup\{k'\}\in\cE_1$. But $P=(P'\backslash\{k\})\cup\{k'\}$ is in $\cE$ and $k\in P'$. So, $k''=k$ and $P''=P$. As a consequence, $\overline\cE$ is a fundamental set of type $(M,n)$ which verifies the $SEU$ property with cardinal strictly smaller than $\cE$. Applying the induction hypothesis on $\overline\cE$, we have the  decomposition of $\cE$ we were looking for.  
\end{demo}

\vskip 2mm

\begin{rmq} The decomposition of the previous proposition induces a decomposition of the vertex set of $\Gamma$. In the proof, we showed that an element of  $\cE_j$ is related only to other elements in $\cE_j$. Consequently, each set $\cE_j$ is the vertex set of a connected component of $\Gamma$. This also implies that this decomposition is unique up to order. We call \emph{connected components} of $\cE$ the sets  $\cE_j$.   
\end{rmq}

\vskip 4mm

\begin{cor} \label{min} Let $\cE$ be a fundamental set of type $(M,n)$ and $\Gamma$ its replacement graph. We assume that $\cE$ verifies the $SEU$ property. Then, the following assertions are equivalent:

\begin{enumerate}
	\item $\cE$ is minimal for the $SEU$ property.
	\item $\cE$ has only one connected component.
	\item $\Gamma$ is connected.
\end{enumerate}
\end{cor}

\vskip 2mm

\begin{rmq} Using propositions \ref{pcomplex} and \ref{facettes}, we can describe $\Gamma$ in terms of $\cP$. The vertices of $\Gamma$ correspond to facets of $\cP$ and two vertices are related if and only if they have a common ridge.
\end{rmq}

\vskip 4mm

Then, we recall the following definition:

\vskip 2mm

\begin{defn} Let $K$ be a simplicial complex. $K$ is a \emph{pseudo-manifold} if the two following properties are fulfilled:
\begin{enumerate}
	\item every ridge of $K$ is contained in exactly two facets.
	\item for all facets $\sigma,\tau$ of $K$, there exists a chains of facets $\sigma=\sigma_1,$ $\dots,$ $\sigma_n=\tau$ of $K$ such that $\sigma_i\cap\sigma_{i+1}$ is a ridge of $K$ for every $i\in\{0,\dots,n-1\}$.
\end{enumerate}
\end{defn}

\vskip 2mm

For instance, every simplicial sphere is a pseudo-manifold. More generally, a triangulation of a manifold (that is, a simplicial complex whose realization is homeomorphic to a topological manifold) is also a pseudo-manifold. Now, the proposition below shows that the notion of pseudo-manifold is exactly the combinatorial property of $\cP$ which characterizes the fact that $\cE$ is minimal for the $SEU$  property:

\vskip 2mm

\begin{prop} \label{pseudomanifold} Let $\cE$ be a fundamental set with $n>M$. Then, $\cP$ is a pseudo-manifold if and only if $\cE$ is minimal for the $SEU$ property.
\end{prop}

\begin{demo} Let assume that $\cE$ is minimal for the $SEU$ property. This implies that every ridge of $\cP$ is contained in exactly two facets  (cf. proposition \ref{crete}). Now, let $\sigma,\tau$ be two distinct facets of $\cP$. So, $P=\sigma^c$ and $Q=\tau^c$ are two fundamental subsets. By minimality for the $SEU$ property, $\Gamma$ is connected (cf. corollary \ref{min}). Consequently, there exists a sequence $P_0=P,P_1,\dots ,P_r=Q$ of fundamental subsets such that $P_i$ and $P_{i-1}$ differ by exactly one element. We denote $R_i$ the acceptable subset $P_{i-1}\cup P_i$ with $M+1$ elements. Its complement $R_i^c$ is thus a face of $\cP$ with $n-M-1=dim(\cP)$ elements. If we put $\sigma_i=P_i^c$, we have $R_i^c=\sigma_{i-1}\cap\sigma_i$ so $\sigma_0=\sigma,\dots ,\sigma_r=\tau$. Consequently, $\cP$ is a pseudo-manifold.
\\Conversely, we assume that $\cP$ is a pseudo-manifold. Then, thanks to proposition \ref{crete}, $\cE$ verifies the $SEU$ property. Moreover, $\cE$ will be minimal for this property if and only if $\Gamma$ is connected (cf. corollary \ref{min}). Let $\sigma,\tau$ be two distinct elements of $\cE$. Then $\sigma^c$ and $\tau^c$ are facets of $\cP$. Since $\cP$ is a pseudo-manifold, there exists a sequence $\sigma_0=\sigma^c$, $\sigma_1,\dots,\sigma_r=\tau^c$ of facets of $\cP$ such that for every $i$, $\sigma_i$ and $\sigma_{i-1}$ share a ridge of $\cP$. This means that $\sigma_i^c$ and $\sigma_{i-1}^c$ are fundamental subsets $\cE$ which differ only by an element, and consequently, $\sigma_i^c$ and $\sigma_{i-1}^c$ are related in $\Gamma$. This implies that $\Gamma$ is connected and $\cE$ is minimal for the $SEU$ property.  
\end{demo}

\vskip 2mm

\begin{rmq} The case where $n=M$ corresponds to $\cP=\{\emptyset\}$. This is not a pseudo-manifold since the only facet is $\emptyset$ and it does not contain any simplex with dimension strictly smaller.
\end{rmq}

\vskip 2mm

Finally, we prove the following proposition which is the motivation for the study of the associated complex:

\vskip 2mm

\begin{prop} Let $(\cE,l)$ be a good system associated to a LVM manifold and ${\mathfrak P}$ its associated polytope. Then the associated complex $\cP$ of $\cE$ is combinatorially equivalent to the boundary of the dual of ${\mathfrak P}$.
\end{prop}

\begin{demo} Since $(\cE,l)$ is a good system associated to a LVM manifold, there exists a bounded component $O$ in $\mathcal O$ such that $\cE$ is exactly the set of subsets $Q$ of $\{1,\dots,n\}$ with $(2m+1)$ for cardinality such that $O$ is included in the convex hull of $(l_q,q\in Q)$. Up to a translation, (whose effect on the action is just to introduce an automorphism of $\C^m\times\C^*$ and so does not change the action, cf. \cite{B}), we can assume that $\cE$ is exactly the set of subsets $P$ of $\{1,\dots,n\}$ with $(2m+1)$ for cardinality such that the convex hull of $(l_p,p\in P)$ contains $0$. In this setting, according to the formula $(7)$ on the page $65$ of \cite{BM}, the boundary of ${\mathfrak P}$ is combinatorially characterized as the set of subsets $I$ of $\{1,\dots,n\}$ verifying 

\[
I \in {\mathfrak P} \Leftrightarrow 0\in Conv(l_k,k\in I^c)
\]

\vskip 2mm

So, $I$ is a subset of  ${\mathfrak P}$ if and only if $I^c$ is acceptable, i.e. $I\in \cP$. As a consequence, from a viewpoint of set theory, ${\mathfrak P}$ and $\cP$ are the same set. We claim that the orders on these sets ${\mathfrak P}$ and $\cP$ are reversed. On the one hand, the order on $\cP$ is the usual inclusion (as for every simplicial complex). On the other hand, we recall the order on the face poset of ${\mathfrak P}$ given in \cite{BM}: every $j$-face of $P$ is represented by a $(n-2m-1-j)$-tuple. So, facets of ${\mathfrak P}$ are represented by a singleton and vertices by a $(n-2m-2)$-tuple. Moreover, a face represented by $I$ is contained in another face represented by $J$ if and only if $I\supset J$. So, combinatorially speaking, the poset of ${\mathfrak P}$ is $(\cP,\supset)$, which prove the claim. Finally, the poset for the dual ${\mathfrak P}^*$ is $(\cP,\subset)$, and the proof is completed.  
\end{demo}


\section{Condition $(K)$}

The previous section shows that $\cP$ is exactly the object we are looking for to generalize the associated polytope of a LVM manifold to the case of LVMB manifolds. We now study its properties. Our first main goal is to prove the following theorem:

\vskip 2mm

\begin{thm} Let $(\cE,l)$ be a good system of type $(2m+1,n)$. Then $\cP$ is a simplicial $(n-2m-2)$-sphere.
\end{thm}

\vskip 2mm

\begin{rmq} The theorem is trivial in the LVM case since the associated complex $\cP$ is a polytope.
\end{rmq}

\vskip 4mm

To prove the previous theorem, we have to focus on good systems which verify an additional condition, called condition $(K)$: 

\vskip 2mm

\begin{center}
$(K)$ \hskip 2mm There exists a real affine automorphism  $\phi$ of $\C^m=\R^{2m}$ such that $\lambda_j=\phi(l_j)$ has coordinates in $\Z^{2m}$ for every $j$. 
\end{center}

\vskip 4mm

For instance, if all coordinates of $l_i$ are rational, then $(\cE,l)$ verifies condition (K). Note that the imbrication condition is an open condition. As a consequence, it is sufficient to prove the previous theorem for good systems verifying the condition $(K)$. Indeed, since $\Q^n$ is dense in $\R^n$, a good system $(\cE,l)$ which does not verify the condition $(K)$ can be replaced by a good system which verifies the condition, with the same associated complex $\cP$.

\vskip 4mm

The main interest of condition $(K)$ stands in the fact that we can associate to our holomorphic acceptable action an algebraic action (called \emph{algebraic acceptable action}) of $(\C^*)^{2m+1}$ on $\C^n$ (or an action of $(\C^*)^{2m}$ on $\PP^{n-1}$):

\vskip 4mm

Let $(\cE,l)$ be a fundamental set of type $(2m+1,n)$ verifying condition $(K)$. We set $l_j=a_j+ib_j,a_j,b_j\in\Z^m$ for every $j$ and $a_j=(a_j^1,\dots,a_j^m)$. We can define an action of $(\C^*)^{2m+1}$ on $\C^n$ by setting: $\forall u\in\C^*,t,s\in\left(\C^*\right)^m,$ $z\in{\C}^n$, we put

\[
\left(u,t,s\right)\cdot z=\left(u \ t_1^{a_1^1}\dots t_m^{a_1^m}\ s_1^{b_1^1}\dots s_m^{b_1^m}\ z_1,\dots, u\ t_1^{a_n^1}\dots t_m^{a_n^m}\  s_1^{b_n^1}\dots s_m^{b_n^m}\ z_n\right)
\]

\vskip 2mm

Using the notation $X_{2m+1}^{\widetilde{l_j}}$ for the character of $(\C^*)^{2m+1}$ defined by $\widetilde{l_j}=(1,l_j)$, we can restate the formula describing the acceptable algebraic action by:

\vskip 2mm

\[
t \cdot z=\left(X_{2m+1}^{\widetilde{l_1}}(t) z_1,\dots, X_{2m+1}^{\widetilde{l_n}}(t)z_n \right)\ \forall t\in (\C^*)^{2m+1},\ z\in\C^n
\]

\vskip 4mm 

It is clear that the open set $\cS$ introduced p.\pageref{ouvertS} is invariant by this action. So we can define $X$ as the topological orbit space of $\cS$ by the algebraic action. As for the holomorphic acceptable action, we can define an action of $\left(\C^*\right)^{2m}$ on $\cV$ whose quotient is also $X$. In \cite{CFZ}, proposition 2.3, it is shown that the holomorphic acceptable action of $\C^m$ on $\cV$ can be seen as the restriction of the algebraic acceptable action to a closed cocompact subgroup $H$ of $(\C^*)^{2m}$. As a consequence, we can define an action of $K=(\C^*)^{2m}/H$ on $\cN$ whose quotient can be homeomorphically identified with $X$.

\vskip 4mm

The principal consequence of this result is the following:

\vskip 2mm

\begin{prop} $X$ is Hausdorff and compact.
\end{prop}

\begin{demo} Let $p:{\cN}\rightarrow X$ be the canonical surjection. $K$ is a compact Lie group so $p$ is a closed map (cf. \cite{Br}, p.38). Consequently, $X$ is Hausdorff. Finally, since $p$ is continuous and $\cN$ is compact, we can conclude that $X$ is compact.  
\end{demo}

\vskip 4mm

Another important consequence for the sequel of the article is that the algebraic action on $\cS$ (or $\cV$) is closed. Moreover, since every complex compact commutative Lie group is a complex compact torus (i.e. a complex Lie group whose underlying topological space is $(S^1)^n$, cf. \cite{L}, Theorem 1.19), we see that $K$ is a complex compact torus.

\vskip 4mm

Using an argument of \cite{BBS}, we show that:

\vskip 2mm

\begin{prop} $t\in(\C^*)^{2m}$ is in the stabilizer of $[z]$ if and only if $\forall i,j\in I_z$, we have 

\[
X_{2m+1}^{l_i}(t)=X_{2m+1}^{l_j}(t)
\]

\end{prop}

\begin{demo} Let us fix some linear order on $\Z^{2m}$. Up to a permutation of the homogeneous coordinates of $\PP^{n-1}$, we can assume that  $l_j\leq l_{j+1}$ for every $j\in\{1,\dots ,n-1\}$ (notice that such a permutation is an equivariant automorphism of $\PP^{n-1}$). We set $j_0=min(I_z)$ the smallest index of nonzero coordinates of $z$. Then, for every $t$ which stabilizes $[z]$, we have

\[
\left[z\right]=t\cdot\left[z\right]=\left[X_{2m+1}^{l_j}(t)z_j\right]
\]

so

\[
\left[z\right]=\left[0,\dots,0,X_{2m+1}^{l_{j_0}}(t)z_{j_0},\dots ,X_{2m+1}^{l_n}(t)z_n\right]
\]

\vskip 2mm 

Consequently, we have 
 
\[
\left[z\right]=\left[0,\dots ,0,z_{j_0},\dots ,X_{2m+1}^{l_n-l_{j_0}}(t)z_j\right]
\]

\vskip 2mm

In particular, we have $X_{2m+1}^{l_j-l_{j_0}}(t)z_j=z_j$ for every $j$. If $j\in I_z$, then  $X_{2m+1}^{l_j-l_{j_0}}(t)=1$.  
\end{demo}

\vskip 2mm

\begin{rmq} In \cite{BBS}, it is shown that $\left[z\right]$ is a fixed point for the algebraic acceptable action if and only if $\forall i,j\in I_z$, we have $l_i=l_j$.
\end{rmq}

\vskip 4mm

Consequently, we have:

\vskip 2mm

\begin{prop}\label{stab} Every element of $\cV$ has a finite stabilizer for the algebraic acceptable action.
\end{prop}

\begin{demo} First, we recall that an element $[z]$ of $\cV$ has at least $(2m+1)$ nonzero coordinates. The index set $I_z$ of these coordinates is an acceptable subset and, by definition, contains a fundamental subset $P$. As a consequence, since $(\cE,l)$ is supposed to be an acceptable system, the set $(l_p)_{p\in P}$ spans $\C^m$ as a real affine space. Up to a permutation, we can assume that $P=(1,2,\dots ,2m+1)$. In this case, we have $X_{2m+1}^{l_j-l_{2m+1}}(t)=1$ for every $j=1,\dots ,2m$ and every $t$ in the stabilizer of $[z]$. We put $L_j=l_j-l_{2m+1}$. Writing $t_j=r_je^{2i\pi\theta_j}$, we get that $r=(r_1,\dots ,r_{2m+1})$ and $\theta=(\theta_1,\dots ,\theta_{2m+1})$ verify the following systems: 

\vskip 2mm

\[
M.ln(r)=0,\ M.\theta\equiv0\ [1]
\]

\vskip 2mm

where $M=(m_{i,j})$ is the matrix defined by $m_{i,j}=L_j^i$ and $ln(r)=(ln(r_1),\dots,ln(r_{2m+1}))$. The system is acceptable, so $(l_1,\dots ,l_{2m+1})$ spans $\C^m$ as a real affine space, which means exactly that the real matrix M is invertible. 

\vskip 4mm

Consequently, $ln(r)=0$ (i.e $|t_j|=1$ for all j) and $\theta\equiv0\ [det(M^{-1})]$. So, $r_je^{2i\pi\theta_j}$ can take only a finite number of values. As a conclusion, the stabilizer of $[z]$ is finite, as claimed.  
\end{demo}

\vskip 2mm

\begin{rmq} An analogous proof shows that the stabilizer of $z\in \cS$ is finite, too. 
\end{rmq}

\vskip 2mm
\section{Connection with toric varieties}

In this section, our main task is to recall that $X$ is a toric variety and to compute its fan. When this fan $\Sigma$ is simplicial, we can construct a simplicial complex $K_\Sigma$ as follows: we denote $\Sigma(1)$ the set of rays of $\Sigma$ and order its elements by $x_1,\dots,x_n$. Then, the complex $K_\Sigma$ is the simplicial complex on $\{1,\cdots,n\}$ defined by:

\[
\forall\ J\subset\{1,\cdots,n\},\hskip 3mm (\ J\in K_\Sigma \Leftrightarrow\ pos(x_j,j\in J)\in \Sigma\ ) 
\]    

\vskip 4mm

This complex $K_\Sigma$ is \emph{the underlying complex} of $\Sigma$. We recall a theorem which will be very important in the sequel:

\vskip 2mm

\begin{prop} \label{sphere} Let X be a normal separated toric variety and $\Sigma$ its fan. We suppose that $\Sigma$ is simplicial. Then, the three following assertions are equivalent:
\begin{enumerate}
\item $X$ is compact.
\item $\Sigma$ is complete in $\R^n$.
\item The simplicial complex underlying $\Sigma$ is a $(n-1)$-sphere.
\end{enumerate}
\end{prop}

\vskip 2mm

\textbf{Convention.}$-$ In the sequel, we will assume that a toric variety is separated and normal.

\subsection{Toric varieties}

To begin with, it is clear that $\cS$ and $\cV$ are   toric varieties. As explained in \cite{CLS}, one associates to a   toric variety with lattice of one-parameter subgroups $N$, a fan $\Sigma$ in the real vector space $N_\R=N\otimes\R$ whose cones are rational with respect to the lattice $N$\footnote{a cone $\sigma$ is said to be rational for $N$ if there a family $S$ of elements of $N$ such that $\sigma=pos(S)$. Moreover, if $S$ is free in $N_\R$, we say that $\sigma$ is simplicial}. So, we can compute the fan associated to $\cS$:

\vskip 2mm

\begin{prop} \label{eventail} Let $(e_i)_{i=1}^n$ be the canonical basis for $\R^n$. Then, the fan describing $\cS$ in $\R^n$ is 

\[
\Sigma(\cS)=\{\ pos(e_i,i\in I)/\ I \in \cP\}
\]

\end{prop}

\vskip 2mm

\begin{demo} To compute the fan of a toric variety, one has to calculate limits for its one-parameter subgroups. The embedding of $(\C^*)^n$ in $\cS$ is the inclusion and the one-parameter subgroups of $(\C^*)^n$ have the form  $\lambda^m(t)=(t^{m_1},\dots,t^{m_n})$ with $m=(m_1,\dots,m_n)\in\Z^n$.
\\So, the limit of $\lambda^m(t)$ when t tends to $0$ exists in $\C^n$ if and only if $m_i\geq 0$ for every $i$. In this case, the limit is $\epsilon=(\epsilon_1,\dots,\epsilon_n)$ with $\epsilon_j=\delta_{m_j,0}$ (Kronecker's symbol).
\\Of course, the limit has to be in $\cS$, which implies for $I_\epsilon$ to be acceptable. But $I_\epsilon=\{j/m_j=0\}$ so the condition means exactly that $\{j/m_j>0\}$ belongs to $\cP$.  
\end{demo}

\vskip 2mm

\begin{exempl} For the fundamental set 

\[
\cE=\left\{\{1,2,5\},\{1,4,5\},\{2,3,5\},\{3,4,5\}\right\}
\] 

\vskip 2mm

of example \ref{exemple}, we have 

\[
{\cS}=\left\{z/(z_1,z_3)\neq 0, (z_2,z_4)\neq 0, z_5\neq 0 \right\}
\]

\vskip 2mm

As a consequence, the fan $\Sigma(\cS)$ is the fan in $\R^5$ whose facets are the $2$-dimensional cones $pos(e_1,e_2),$ $pos(e_1,e_4),$ $pos(e_2,e_3)$ and $pos(e_3,e_4)$ (where $e_1,\dots,e_5$ is the canonical basis of $\C^5$).  
\end{exempl}

\vskip 2mm

\begin{rmq}
We can also easily compute orbits of $\cS$ for the action of $(\C^*)^n$. For $I\subset\{1,\dots,n\}$, we set $O_I=\{z\in\C^n/I_z=I^c\}$. Then if $z\in\cS$, its orbit is $O_{I_z^c}$.
\\So, we obtain the partition: $\cS=\displaystyle{\bigsqcup_{I\in\cP} O_I}$

\vskip 2mm

Moreover, in the orbit-cone correspondence between $\cS$ and $\Sigma(\cS)$ (cf. \cite{CLS} ch.3), $O_{I}$ corresponds to $\sigma_I=pos(e_j/j\in I)$.
\end{rmq}

\vskip 2mm

\begin{rmq} In quite the same way, we can show that the fan of $\cV$ is $\Sigma(\cV)=\{pos(e_i,i\in I)/I\in \cP\}$ in $\R^{n-1}$, with $(e_1,\dots,e_{n-1})$ defined as the canonical basis of $\R^{n-1}$ and $e_n=-(e_1+\dots+e_{n-1})$.
\end{rmq}

\vskip 2mm

In \cite{CFZ}, it is proven that the quotient $X$ of the acceptable algebraic action of $\left(\C^*\right)^{2m+1}$ on $\cS$ is a compact toric variety. In the next section, we will detail the construction with in order to identify its group of one-parameter subgroups and the structure of its fan.

\vskip 4mm

\begin{exempl} For the good system $(\cE,l)$ with 
\[
\cE=\left\{\{1,2,5\},\{1,4,5\},\{2,3,5\}\{3,4,5\}\right\}
\] 

and $l_1=l_3=1$, $l_2=l_4=i$ and $l_5=0$, the algebraic action is 

\[
\left(\alpha,t,s\right)\cdot z=\left(\alpha t z_1,\alpha s z_2,\alpha t z_3,\alpha s z_4,\alpha z_5\right)
\]

\vskip 2mm

Using the automorphism of $\left(\C^*\right)^3$ defined by $\phi(\alpha,t,s)=(\alpha,\alpha t,\alpha s)$, we can see that the quotient $X$ of the algebraic action is also the quotient of $\cS$ by the action defined by 

\[
(\alpha,t,s)\cdot z=(t z_1,s z_2,t z_3,s z_4,\alpha z_5)
\]

\vskip 2mm

so $X$ is the product $\PP^1\times\PP^1$.
\end{exempl}

\vskip 4mm

To conclude this section, let $f:(\C^*)^{2m+1}\longrightarrow (\C^*)^n$ be the map defined by

\vskip 2mm

\[
f(u,t,s)=(X_{2m+1}^{\widetilde{l_1}}(u,t,s),\cdots,X_{2m+1}^{\widetilde{l_n}}(u,t,s))
\]

\vskip 2mm

The algebraic acceptable action on $\C^n$ is just the restriction to $Im(f)$ of the natural action of the torus $(\C^*)^n$ on $\C^n$. 

\vskip 2mm

\begin{prop} $Ker(f)$ is finite.
\end{prop}

\vskip 2mm

\begin{demo} Let $z$ be the point $(1,\dots,1)$. Then, $I_z=\{1,\dots,n\}$ and $z\in\cS$. The stabilizer of $z$ for the action of $\left(\C^*\right)^{2m+1}$ is exactly $Ker(f)$ then the remark following the proposition \ref{stab} implies that $Ker(f)$ is finite.
\end{demo}

\vskip 2mm

\begin{rmq} Generally, $f$ is not injective. For instance, if we consider $\cE=\{\{1,2,4\},\{2,3,4\}\}$ and $l_1=1,l_2=i,l_3=p,l_4=-1-i$, where $p$  is a nonzero positive integer. $(\cE,l)$ is a good system and for $p=4$, $f$ is not injective. Note that for $p=3$, $f$ is injective.
\end{rmq}

\vskip 2mm

\begin{defn} We define $T_N$ as the quotient group $\left(\C^*\right)^n/Im(f)$
\end{defn}

\vskip 2mm

We recall that $\left(\C^*\right)^n$ is included in $\cS$ as a Zariski open subset. $(\C^*)^n$ is invariant by the action of $(\C^*)^{2m+1}$. This implies that $T_N$ can be embedded in $X$ as an dense open subset. Moreover, the action of $(\C^*)^n$ on $\cS$ commutes with the action of $(\C^*)^{2m+1}$, so the action of $T_N$ on itself can be extended to an action on $X$.

\subsection{The algebraic torus $T_N$}

We denote $F:\C^{2m+1}\rightarrow\C^n$ the linear map defined by 

\[
F\left(U,T,S\right)=\left(U+<a_j,T>+<b_j,S>\right)_j
\]

\vskip 4mm

The matrix of $F$ has $(1,l_j)$ as $j$-th row so $F$ has maximal rank. So, $F$ is injective. The family $f_j=F(e_j)$, with $(e_1,\dots ,e_{2m+1})$ the canonical basis of $\C^{2m+1}$ is a basis for $Im(F)$. We notice that each $f_j$ has integer coordinates. We complete this basis to a basis $(f_1,\dots ,f_n)$ of $\C^n$ with integer coordinates. Next, we define the map $G:\C^n\rightarrow\C^{n-2m-1}$ by linearity and 
$G(f_j)=\left\{\begin{array}{cl}
										0 	 & j\in\{1,\dots,2m+1\}\\ 
										g_j  & otherwise\\ 
							\end{array}
				\right.$ 
				
\vskip 2mm				

(with $(g_{2m+2},\dots ,g_n)$ the canonical basis of $\C^{n-2m-1}$). 

\vskip 2mm

It is clear by construction that the following sequence is exact:

\[ 
\xymatrix{
				0 \ar[r] & \C^{2m+1} \ar[r]^F & \C^n \ar[r]^G & \C^{n-2m-1} \ar[r] & 0
   			}
\]

Moreover, we have 

\[
F(\Z^{2m+1})\subset \Z^n,\ G(\Z^n)\subset \Z^{n-2m-1}
\]

\vskip 2mm

Let $t\in\left(\C^*\right)^{2m+1}$ and $T$ be some element of $\C^{2m+1}$ such that $t=exp\left(T\right)$. We put $g\left(t\right)=exp\left(G\left(T\right)\right)$. The previous remark has for consequence that $g$ is well defined. Moreover, we have:

\vskip 2mm

\begin{prop} $g$ is a group homomorphism and the following diagram is commutative:

\[
\xymatrix{
    \left(\C^*\right)^{2m+1} \ar[r]^-{f}  & \left(\C^*\right)^n \ar[r]^-{g}  & \left(\C^*\right)^{n-2m-1} \\
    \C^{2m+1} \ar@{->>}[u]_-{exp}  \ar@{^{(}->}[r]^-{F} &  \C^n \ar@{->>}[u]_-{exp} \ar@{->>}[r]^-{G} & \C^{n-2m-1} \ar@{->>}[u]_-{exp}  
  }
\]

\end{prop}

\vskip 2mm

Finally, we obtain:

\vskip 2mm

\begin{prop} $g$ is surjective and $Ker(g)=Im(f).$
\end{prop}

\begin{demo} The surjectivity of $g$ is clear (since $G$ and $exp$ are surjective).
\\So, we have only to show that $Ker(g)=Im(f)$.
By the construction of $F$ and $G$ and by commutativity of the previous diagram, we have for every $t=exp(T)$, $g\circ f(t)=g\circ f\circ exp(t)=exp\circ G\circ F(T)=exp(0)=1$ so $Im(f)\subset Ker(g)$. Conversely, let $t$ belong to $Ker(g).$ We put $t=exp(T)$, for some $T\in \C^n$ and we have $1=g(t)$ so $exp(G(T))=1$. As a consequence, $G(T)\in 2i\pi\Z^{n-2m-1}$, i.e. 

\[
G(T)=(2i\pi q_{2m+2},\dots ,2i\pi q_n)
\]

and

\[
G(T)=2i\pi q_{2m+2}g_{2m+2}+\dots 2i\pi q_ng_n=G(2i\pi q_{2m+2}f_{2m+2}+\dots 2i\pi q_nf_n)
\]

\vskip 2mm

So $T-(2i\pi q_{2m+2}f_{2m+2}+\dots 2i\pi q_nf_n)\in Ker(G)=Im(F)$. We have that 

\[
T=\lambda_1f_1+\dots\lambda_{2m+1}f_{2m+1}+2i\pi q_{2m+2}f_{2m+2}+\dots 2i\pi q_nf_n
\]

\vskip 2mm

Finally, $T=F(\lambda_1e_1+\dots\lambda_{2m+1}e_{2m+1})+2i\pi(q_{2m+2}f_{2m+2}+\dots q_nf_n)$, which implies that 

\[
t=exp(F(\lambda_1e_1+\dots\lambda_{2m+1}e_{2m+1}))
\]

\vskip 2mm

which means that

\[
t=f(exp(\lambda_1e_1+\dots\lambda_{2m+1}e_{2m+1}))\in Im(f)\
\]

\end{demo}

\vskip 3mm

In particular, $T_N=(\C^*)^n/Im(f)$ is isomorphic to $(\C^*)^{n-2m-1}$ (and, as claimed in the previous section, $X$ is a toric variety). We denote $\overline{g}$ for the isomorphism between $T_N$ and $(\C^*)^{n-2m-1}$ induced by $g$.
\vskip 2mm

\begin{defn} We denote $\lambda_T^u$ the one-parameter subgroup of $T_N$ defined by $\lambda_T^u=\overline{g}^{-1}\circ \lambda_{n-2m-1}^u$. Since  $\overline{g}$ is an isomorphism, every one-parameter subgroup of $T_N$ has this form.
\end{defn}

\vskip 2mm

\textbf{Notation:} Let $\phi:T_1\rightarrow T_2$ be a group homomorphism between two algebraic tori $T_1$ and $T_2$. We will denote $\phi^*$ the morphism between the groups of one-parameter subgroups of $T_1$ and $T_2$ induced by $\phi$: $\phi^*(\lambda)=\phi\circ\lambda$.

\vskip 2mm

The group of one-parameter subgroups of $\left(\C^*\right)^n$ is $\{\lambda_n^u/u\in\Z^n\}$ which we will identify with $\Z^n$ via $u\leftrightarrow \lambda_n^u$. Via this map, the map $F$ and $G$ are exactly the morphisms induced by $f$ and $g$ (respectively):

\vskip 2mm

\begin{prop} With the above identification, we have $F=f^*$ and $G=g^*$. Precisely: 
\begin{enumerate} 
\item For every $v$ in $\Z^{2m+1}$, $f\circ\lambda_{2m+1}^v$ is the one-parameter subgroup $\lambda_n^{F(v)}$ of $\left(\C^*\right)^n$.
\item For every $v$ in $\Z^n$, $g\circ\lambda_n^v$ is the one-parameter subgroup $\lambda_{n-2m-1}^{G(v)}$ of $\left(\C^*\right)^{n-2m-1}$.
\end{enumerate}
\end{prop}

\begin{demo} Let us take some $t\in\C^*$. 
\\1) We have $\lambda_{2m+1}^v(t)=(t^{v_1},...,t^{v_{2m+1}})$ so \[f\circ\lambda_{2m+1}^v(t)=(t^{v_1+v_2a_j^1+\dots+v_{m+1}a_j^m+\dots+v_{2m+1}a_j^m})_j=(t^{<v,\widetilde{l_j}>})_j=\lambda_n^{F(v)}(t)\]

2) Let us pick some $T\in\C^n$ such that $t=exp(T)$. We put $w=G(v)=(w_1,...,w_{n-2m-1})$. If we denote $g_j$ (resp. $G_j$) for the coordinate functions of $g$ (resp. $G$) in the canonical bases, we can easily verify that $g_j\circ exp=exp\circ G_j$ and that $w_j=G_j(v)$ for every $j$.

\vskip 2mm

Next, we have $\lambda^v_n(t)=(e^{Tv_1},$ $\dots,$ $e^{Tv_n})$, so

\[
g\circ\lambda^v_n(t)=(e^{G_1(Tv)},\dots,e^{G_{n-2m-1}(Tv)})
\]

and $g\circ\lambda^v_n(t)=\lambda_{n-2m-1}^w(t)$.   
\end{demo}

\vskip 2mm

We would like to identify the group $N$ of one-parameter subgroups of $T_N$ with some lattice in $\C^n/Im(F)$. A natural candidate is described as follows: Let $\Pi$ be the canonical surjection $\C^n\rightarrow \C^n/Im(F)$ and $\overline{G}:\C^n/Im(F)\rightarrow\C^{n-2m-1}$ be the linear isomorphism induced by $G$. Notice that $\overline{G}$ is also a $\Z$-module isomorphism between $\Z^{n-2m-1}$ and $\Pi(\Z^n)$. If $u$ belongs to $\Z^n$, we set $\lambda_T^{\Pi(u)}$ for the one-parameter subgroup $\lambda_T^{G(u)}$ (notice that this definition makes sense since $Im(F)=Ker(G)$). As a consequence, we have: $N=\{\lambda_T^{\Pi(u)}/u\in\Z^n\}$ which can be identified with $\Pi(\Z^n)=\Z^n/Im(F).$ 

\vskip 3mm 

Now, we can define an "exponential" map between $\C^n/Im(F)$ and $T_N$: we define $exp:\C^n/Im(F)\rightarrow T_N$ by setting $exp(\Pi(z))=\pi\circ(exp(z))$ for every element $z$ in $\C^n$. The fact that $Ker(g)=Im(f)$ implies that this map is well defined.(Alternatively, we can define this exponential as the map $\overline{g}^{-1}\circ exp\circ\overline{G}$). By construction, we have $\pi\circ exp=exp\circ\Pi$ and $exp\circ\overline{G}=\overline{g}\circ exp$. Moreover, for every $v\in\Z^n$, $\pi\circ\lambda_n^v=\lambda_T^{\Pi(v)}$, that is, that $\Pi=\pi^*$.

\vskip 2mm
%
%

To sum up, we have the following commutative diagram:
\label{suitesexactes}
\[
\xymatrix{
		0 \ar[r] & \C^{2m+1} \ar@{^{(}->}[rr]^F \ar[rd]^F \ar@{->>}[dd]^{exp} & & \C^n \ar@{->>}@/_2pc/[ll]^{\widetilde{F}} \ar@{->>}[rr]^G \ar@{->>}[rd]^\Pi \ar@{->>}[dd]^{exp} & & \C^{n-2m-1} \ar[r] \ar@{->>}[dd]^{exp} \ar@{^{(}->}@/_2pc/[ll]^{\widetilde{G}} & 0\\		
   	& & Im(F) \ar@{^{(}->}[ru]^i \ar@{->>}[dd] & & \C^n/Im(F) \ar[ru]_{\overline{G}} \ar@{->>}[dd]^{exp}  \\
    & (\C^*)^{2m+1} \ar@{-}[r]^f \ar@{->>}[rd]^f & \ar[r] & (\C^*)^n \ar@{-}[r] \ar@{->>}[rd]^\pi & \ar@{->>}[r]^-{g} & (\C^*)^{n-2m-1} \ar[r] & 1\\
    & & Im(f) \ar@{^{(}->}[ru]^i & & T_N \ar[ru]_{\overline{g}} \ar@{=}[d] \\
    & & & & (\C^*)^n/Im(f)
}
\]


\subsection{Study of the projection $\pi$}

In this last section, we are in position to prove the first part of theorem \ref{thmprinc}. Let $\Sigma$ denote the fan in $N_\R=N\otimes\R$ associated to $X$ (this fan exists since $X$ is separated and normal). In order to use proposition \ref{sphere}, we have to prove that $\Sigma$ is simplicial. According to \cite{CLS}, we just have to prove that $X$ is an orbifold. 

\vskip 4mm

As shown previously, the holomorphic acceptable action on $\cS$ is free and the algebraic action on $\cS$ has only finite stabilizers. Consequently, every stabilizer for the action of $K$ on $\cN$ is finite. So, $X$ is the quotient of the compact variety $\cN$ by the action of a compact Lie group $K$ and every stabilizer for this action is finite. 

\vskip 4mm

Finally,  we can claim that the map 

\[
\begin{array}{rccl} \phi_n: & K & \longrightarrow  &{\cN} \\ & h & \mapsto & h\cdot n  \end{array}
\]

is proper for every element $n$ in $\cN$ since it is a continuous map defined on a compact set. By Holmann's theorem (cf. \cite{O}, \S 5.1), we get that $X$ is indeed an orbifold. Taking advantage of proposition \ref{sphere}, we have proved that $K_\Sigma$, the underlying complex of the fan $\Sigma$, is a $(n-2m-2)$-sphere.

\vskip 4mm

Now, the following proposition will complete the proof of the theorem \ref{thmprinc}:

\vskip 2mm

\begin{thm} $K_\Sigma$ and $\cP$ are isomorphic simplicial complexes.
\end{thm}

\begin{demo} First of all, $\pi$ is by construction a toric morphism. According to \cite{CLS}, this implies that $\pi$ is equivariant (for the toric actions of $\left(\C^*\right)^n$ and $T_N$), so $\pi$ sends an orbit $O_I$ in $\cS$ to an orbit in $X$. We set $\widetilde{O_I}$ to be the unique orbit in $X$ containing $\pi(O_I)$. Moreover, the map $\pi^*$ (identified with $\Pi$) from $\R^n$ into $N_\R$ preserves the cones.

\vskip 4mm

We can easily show that $\cS$ and $X$ have exactly the same number of orbits, i.e that the quotient of $X$ by its torus $T_N$ is in bijection with the quotient of $\cS$ by its torus $\left(\C^*\right)^n$. Since $\pi$ is surjective, we get that the assignment  $O_I\rightarrow\widetilde{O}_I$ (induced by $\pi$) is bijective. As a consequence, $\Pi$ induces a bijection between the cones of $\Sigma(\cS)$ and those of $\Sigma$

\vskip 4mm

If $\sigma$ is a cone belonging to $\Sigma(\cS)$, we denote $O(\sigma)$ the orbit in $\cS$ associated to $\sigma$ (cf \cite{CLS}, ch.3). Particularly, $O(\sigma_I)=O_I$. We will also denote in the same way the orbits in $X$. Moreover, the image of the cone $\sigma$ by $\Pi$ will be denoted $\widetilde{\sigma}$. So, we have $\widetilde{O(\sigma)}=O(\widetilde{\sigma})$.
\\Still from \cite{CLS}, $\pi$ preserves the partial order of faces: if $\tau$ is a face of $\sigma$ in $\Sigma(\cS)$, then $\widetilde{\tau}$ is a face of $\widetilde{\sigma}$. 

\vskip 4mm 

On the other hand, a slight modification of the proof of the fact that $(\C^*)^n/Im(f)$ (i.e. $\pi\left(O_\emptyset\right)$) is isomorphic to $\left(\C^*\right)^{n-2m-1}$ shows that $\pi(O(\sigma))$ is isomorphic to $(\C^*)^{n-2m-1-dim(\sigma)}$. At the level of cones, this means that cones of $\Sigma(\cS)$ with the same dimension are sent to cones of $\Sigma$ with the same dimension. In particular, $\Pi$ sends rays to rays. This last property means that $\Pi$ induces a bijection between vertices of $\cP$ and vertices of $K_\Sigma$, bijection which we also denote $\Pi$.
\\It is clear by what precedes that this very last map is an isomorphism of simplicial complexes.  
\end{demo}

\vskip 4mm

As a consequence, $\cP$ is indeed a simplicial sphere. More precisely, what we have is a particular type of simplicial spheres, namely spheres which are the underlying simplicial complex of some complete fan:

\vskip 4mm

\begin{defn} Let $K$ be a $d$-sphere. $K$ is said to be rationally starshaped if there exists a lattice $N$ in $\R^{d+1}$, a point $p_0\in N$ and a realization\footnote{The fact that a $d$-sphere has a realization in $\R^{d+1}$ is an open question (cf. \cite{MW}, \S 5).} $|K|$ for $K$ in $\R^{d+1}$ such that every vertex of $|K|$ belongs to $N$ and every ray emanating from $p_0$ intersects $|K|$ in exactly one point. The realization $|K|$ is said to be \emph{starshaped} and we say that $p_0$ belongs to the \emph{kernel} of $|K|$. 
\end{defn}

\vskip 4mm

\begin{cor} Let $(\cE,l)$ be a good system verifying $(K)$. Then its associated complex $\cP$ is a rationally starshaped sphere.
\end{cor}

\begin{demo} We have already seen that $\cP$ is a simplicial sphere combinatorially equivalent to $K_{\Sigma}$. If we put $\Sigma(1)=\{\rho_1,\dots,\rho_v\}$ for the distinct rays of $\Sigma$ and $u_1,\dots,u_v$ for generators of $\rho_1,\dots,\rho_v$ (respectively) in $N$, then the geometric simplicial complex $C$ whose simplexes are $Conv(u_i,i\in I)$ for $I$ in $\cP$ is obviously a realization of $\cP$ in $\R^{n-2m-1}$ with rational vertices. The point $0$ is in the kernel of $C$ so $\cP$ is rationally starshaped.  
\end{demo}

\vskip 2mm 
\section{Inverse construction}

\subsection{The construction}

In this section, we give a realization theorem: for every rationally starshaped sphere, there exists a good system whose associated sphere is the given one. In \cite{M},\ p.86, the same kind of theorem is proven for simple polytopes and LVM manifolds. One of the main interests of this theorem is that it gives us a clue to answer an open question: does there exist a LVMB manifold which has not the same topology as a LVM manifold? The idea is to use this theorem to construct a LVMB manifold from a rationally starshaped sphere $\cP$ which is not polytopal and that this manifold has a particular topology. For example, we expect that the LVMB manifold coming from the Brückner sphere or the Barnette one (the two $3$-spheres with $8$ vertices which are not polytopal) are good candidates.

\vskip 4mm

Let $\cP$ be a rationally starshaped $d$-sphere with $v$ vertices (up to an isomorphism of simplicial complexes, we will assume that these vertices are $1,2,\dots,v$). So, there exists a lattice $N$ and a realization $|\cP|$ of $\cP$ in $\R^{d+1}$ all of whose vertices belong to $N$. We can assume that $0$ is in the kernel of $|\cP|$ and that $N$ is $\Z^{d+1}$. We denote $x_1,\dots,x_{d+1}$ the vertices of $|\cP|$ corresponding to the vertices $1,2,\dots,v$ of $\cP$ and $p_1,\dots,p_v$ the generators of the rays of $\Z^{d+1}$ passing through $x_1$ (i.e. $p_j$ is the unique generator of the semi-group $\Z^{d+1}\cap[0,x_j)\ $). 

To begin, we suppose that $v$ is even and we put $v=2m$. We denote $\cE$ the set defined by 

\[
\left\{\ P\subset\{1,...,v+d+1\}\middle/\ P^c\ is\ a\ facet\ of\ \cP\ \right\}
\] 

\vskip 4mm

and $\cE_0=\{\{0\}\cup P / P\in \cE\}$. We also denote $A$ the matrix whose columns are $p_1,\dots,p_v$. We label its rows $p^1,\dots,p^{d+1}$, and finally, we have:

\vskip 4mm 

\begin{thm} If $e_1,\dots,e_v$ is the canonical basis of $\R^v$, then 

\[
(\cE_0,(0,e_1,\dots,e_v,-p^1,\dots,-p^{d+1}))
\] 

is a good system of type $(v+1,d+v+2)$ whose associated complex is $\cP$.
\end{thm}

\begin{demo}
First, since $\cP$ is a $d$-sphere, it is a pure complex whose facets have $d+1$ elements. As a consequence, every subset in $\cE_0$ has $v+1$ elements. So, $\cE_0$ is a fundamental set of type $(v+1,v+d+2)$ (as usual, we will denote the set of its acceptable subsets $\cA_0$). Notice that, since every facet of $\cP$ has vertices in $\{1,\dots,v\}$, $\cE_0$ has $d+2$ indispensable elements. Moreover, by definition, the associated complex to $\cE_0$ is clearly $\cP$. Since $\cP$ is a sphere, hence a pseudo-manifold, proposition \ref{pseudomanifold} implies that $\cE_0$ is minimal for the SEU property. We can also notice that the same is true for $\cE$ (with set of acceptable subsets $\cA$).

\vskip 4mm

Secondly, we have to check that the vectors 

\[
0,e_1,e_2,\dots,e_v,-p^1,\dots,-p^{d+1}
\]

"fit" with $\cE_0$ to make a good system. We put $\rho_j=pos(p_j)$ for the ray generated by $p_j,j=1,\dots,v$. We also denote $\Sigma(\cP)$ the fan defined by 

\[
\Sigma(\cP)=\left\{\ pos(p_j,j\in I)\ \middle/\ I\in\cP\ \right\}
\]

\vskip 4mm

Then, by definition, $\Sigma(\cP)$ is a simplicial fan whose underlying complex is $\cP$. As a consequence, $\Sigma(\cP)$ is rational with respect to $N$ and complete (since $\cP$ is a sphere, cf. proposition \ref{sphere}). We set $X$ for the compact  toric variety associated to $\Sigma(\cP)$. Following \cite{Ha}, we will construct $X$ as a quotient of a quasi-affine toric variety by the action of an algebraic torus. 

\vskip 4mm

In what follows, $(e_j)_{j=1,\dots,N}$ is the canonical basis of $\C^N$ (for any $N$). Let $\widetilde{\Sigma}$ be the fan in $\R^{v+d+1}$ whose cones are $pos(e_j,j\in J)$, $J\in\cP$. Obviously, $\widetilde{\Sigma}$ is a non complete simplicial fan whose underlying complex is $\cP$, too. We denote $\widetilde{X}$ the (quasi-affine) toric variety associated to this fan. Then, the open set $\cS=\{z\in \C^{v+d+1}/I_z\in \cA\}$ is exactly the set $\widetilde{X}$. Indeed, the computation of the proof of proposition \ref{eventail} shows that the fan of $\cS$ is the fan in $\R^{v+d+1}$ whose rays are generated by the canonical basis and whose underlying complex is $\cP$. Thus, $\cS$ and $\widetilde{X}$ have exactly the same fan, so they coincide. We also observe that 

\[
\cS_0=\C^*\times \cS=\left\{\ (z_0,z)\in\C^{v+d+2}\ \middle/\ I_{(z_0,z)}\in \cA_0\ \right\}
\] 

\vskip 4mm

has the same fan, but seen in $\R^{v+d+2}$. 

\vskip 2mm

Finally, we define the map $f:(\C^*)^v\rightarrow(\C^*)^{v+d+1}$ by 

\[
f(t)=(t,X_v^{-p^1}(t),\dots,X_v^{-p^{d+1}}(t))
\]

\vskip 4mm

Note that $t=(X_v^{e_1}(t),\dots,X_v^{e_v}(t))$. According to \cite{Ha}, $X$ is the quotient of $\widetilde{X}=\cS$ by the restriction of the toric action of $(\C^*)^{v+d+1}$ restricted to $Im(f)$. It is a geometric quotient because $\widetilde{\Sigma}$ is simplicial.  

\vskip 2mm

Considering $l_0$=$0$,$l_1$=$e_1$,$\dots,l_{2m}$=$e_{2m}$,$l_{v+1}$=$-p^1$,$\dots$,$l_{v+d+1}$=$-p^{d+1}$ as elements of $\C^m$, we can define a holomorphic action of $(\C^*)\times \C^m$ on $\cS_0$ by setting 

\[ 
(\alpha,T)\cdot (z_0,z)=(\alpha e^{<l_j,T>}z_j)_{j=0}^{v+d+1}\hskip 4mm \forall \alpha\in\C^*,t\in \C^m, (z_0,z)\in\cS_0
\] 

\vskip 2mm

It is clear that $(\cE_0,l)$ verifies $(K)$. Moreover, the algebraic acceptable action associated to this system has $X$ for quotient. Indeed, a computation shows that this algebraic action is defined by 

\[
(\alpha,t)\cdot (z_0,z)=(\alpha z_0,f(t)\cdot z)\hskip 4mm  \forall \alpha\in\C^*,t\in(\C^*)^{2m}, (z_0,z)\in\cS_0
\]

\vskip 2mm
 
If we denote $\cV_0=\left\{\ \left[z_0,z\right]\ \middle/\ \left(z_0,z\right)\in\cS_0\ \right\}$ the "projectivization" of $\cS_0$, then the orbit space for the algebraic action is the quotient of $\cV_0$ by the action defined by 

\[
t \cdot \left[z_0,z\right]=\left[z_0,f(t)\cdot z\right]\hskip 4mm  \forall T\in(\C^*)^{2m}, \left[z_0,z\right]\in\cV_0
\]

\vskip 2mm

But $\cV_0=\left\{\ \left[1,z\right]\ \middle/\ z\in \cS\ \right\}$ is homeomorphic to $\cS$ so the claim that $X$ is the orbit space of the action of $\left(\C^*\right)^{2m+1}$ on $\cS_0$ is proved.

\vskip 2mm

The only thing we have to check is that $(\cE_0,l)$ is a good system, that is the orbit space $\cN$ for the holomorphic action is a complex manifold. We have seen that $X$ is a geometric quotient so the action of $\left(\C^*\right)^{2m}$ is proper (see, for example, \cite{BBCMcG}, p.28). As a consequence, the action of $\C^m$ is proper too, and since there are no compact subgroups in $\C^m$ (except $\{0\}$ of course), this action is free. Finally, the action is proper and free so $\cN$ can be endowed with a structure of complex compact manifold.
\end{demo}

\vskip 2mm 

\begin{rmq} 
1) Since $\cN$ is a complex compact manifold, the imbrication condition is fulfilled (cf. Theorem \ref{thmBosio1}). This gives another proof for the fact that the SEU property is also fulfilled.
\\2) The transformation of 

\[
\begin{array}{rcl}
 						 & & (p_1,\dots,p_v,e_1,\dots,e_{d+1})   \\
 \mathrm{to} & & 																		 \\ 
 						 & & (e_1,\dots,e_v,-p^1,\dots,-p^{d+1}) 
\end{array} 						 
\] 

is called a linear transform of $(x_1,\dots,x_v,e_1,\dots,e_{d+1})$ (see \cite{E} for example). In \cite{M}, the construction of a LVM manifold starting from a simple polytope used a special kind of linear transform called Gale (or affine) transform (see. \cite{E})
\end{rmq}

\vskip 2mm 

Now, we suppose that $v=2m+1$ is odd. The construction of  good system whose associated complex is $\cP$ is very similar. However, we have an additional step: we define an action of $(\C^*)^{v+1}$ on $\cS_0$ by \[(t_0,t)\cdot (z_0,z)=(t_0z_0,f(t)\cdot z)=(t_0z_0,X^{e_1}(t)z_1,\dots,X^{-x^{d+1}}(t)z_{v+d+1})\] where $f$ and the action of $Im(f)$ are defined as above, and $e_0$ is the first vector of the canonical basis of $\C^{v+1}$ with coordinates $(z_0,z_1,\dots, z_{v})$. The orbit space for this last action is still $X$. The rest is as above: we define a fundamental set $\cE_*=\{\ \{-1\}\cup E\ /\ E\in \cE_0\ \}$ and we have:

\vskip 2mm 

\begin{thm} If $e_0,\dots,e_v$ is the canonical basis of $\R^{v+1}$, then 

\[
\bigg(\ \cE_*,\ \Big(\ 0,e_0,e_1,\dots,e_{v},(0,-p^1),\dots,(0,-p^{d+1})\ \Big)\ \bigg)
\]

is a good system of type $(v+2,d+v+3)$ whose associated complex is $\cP$.
\end{thm}


\subsection{LVMB manifolds and moment-angle complexes}

\label{cma}

In this section, we follow closely the definitions and notations of \cite{BP}. We use the result of the above section to show that many moment-angle complexes with even dimension can be endowed with a complex structure of a LVMB manifold.

\vskip 2mm

\begin{defn} Let $K$ be a simplicial complex on $\{1,\dots,n\}$ with dimension $d-1$. If $\sigma\subset\{1,\dots,n\}$, we put 

\[
C_\sigma=\left\{\ t\in\left[0,1\right]^n\ \middle/\ t_j=1\ \forall j\notin \sigma\ \right\}
\] 

and 

\[
B_\sigma=\left\{\ z\in\D^n\ \middle/\ \left|z_j\right|=1\ \forall j\notin \sigma\ \right\}
\]

\vskip 2mm

The moment-angle complex associated to $K$ is 

\[
\cZ_{K,n}=\bigcup_{\sigma\in K} B_\sigma
\]

\end{defn}

\vskip 2mm

\begin{exempl} For instance, if $K$ is the boundary of the $n$-simplex, $\widetilde{\cZ}_K$ is the sphere $S^{2n-1}$ (cf. \cite{BP}).
\end{exempl}

\vskip 2mm

In \cite{BP}, Lemma $6.13$, it is shown that if $K$ is a simplicial sphere, then $\cZ_{K,n}$ is a closed manifold. Moreover, let $(\cE,l)$ be a good system with associated complex $\cP$. In \cite{B}, to prove Theorem \ref{thmBosio1} (p.1268 in \cite{B}), Bosio introduces the set $\widehat{M_1'}$ defined by 

\[
\widehat{M_1'}=\left\{\ z\in\D^n \middle/\ J_z \in \cA\ \right\}
\]

\vskip 2mm 

where $J_z=\{\ k\in\{1,\dots,n\}/\ |z_j|=1\ \}$. This set is the quotient of $\cS$ by the restriction of the holomorphic action to $\R_+^*\times\C^m$ and as a consequence, $\cN$ is the quotient of $\widehat{M_1'}$ by the \emph{diagonal action} of $S^1$ defined by:

\vskip 2mm

\[ 
e^{i\theta}\cdot z=(\ e^{i\theta}z_1,\dots,e^{i\theta}z_n\ )\hskip 3mm \forall \theta\in\R,\ z=(z_1,\dots,z_n)\in \widehat{M_1'}\subset \C^n
\]

\vskip 2mm

\begin{prop}\label{complexe} We have $\widehat{M_1'}=\cZ_{\cP,n}$
\end{prop}

\begin{demo} Indeed, we have: 

\[\begin{array}{rcl}
\widehat{M_1'} & = & \left\{\ z\in \D^n\ \middle/\ J_z\in\cA\ \right\} \\
 		    			 & = & \displaystyle{\bigcup_{\tau\in\cA}\left\{\ z\in\D^n\ \middle/\tau\subset J_z\ \right\}} \\
 		    			 & = & \displaystyle{\bigcup_{\sigma\in\cP} \left\{\ z\in\D^n\ \middle/\ \forall j\notin\sigma, \left|z_j\right|=1\ \right\}}  \\
 		    			 & = & \cZ_{\cP,n} \\
\end{array}
\] 

\end{demo}

\vskip 2mm

\begin{prop} \label{complexe2} Let $(\cE,l)$ be a good system with type $(2m+1,n,k)$ and $\cN$ the LVMB associated to this system. If $k>0$, then $\cN$ is homeomorphic to a moment-angle complex.
\end{prop}

\begin{demo} We assume that $n$ is indispensable. Let $\cP$ be the sphere associated to $\cE$ and $\cS$ the open subset of $\C^n$ whose quotient by the holomorphic action is $\cN$. According to proposition \ref{complexe}, the quotient $\widehat{M_1'}=\cS/(\R_+^*\times \C^m)$ can  be identified with $\cZ_{\cP,n}$. 

\vskip 2mm

Let $\phi$ be the map defined by 

\[
\begin{array}{rcccl}
\phi & : & \cZ_{\cP,n} & \rightarrow & \C^{n-1} \\
	   &	 & z					 & \mapsto 		 & \left(\frac{z_1}{z_n},\dots, \frac{z_{n-1}}{z_n} \right) 
\end{array}
\]

Since $n$ is indispensable, we have $|z_n|=1$ for every $\cZ_{\cP,n}$ so $\phi$ is well defined. Moreover, $\phi$ is continuous and a simple calculation shows that $\phi$ is invariant for the diagonal action and $\phi(\cZ_{\cP,n})=\cZ_{\cP,n-1}$. We claim that if $\phi(z)=\phi(w)$, then $z$ and $w$ belong to the same orbit for the diagonal action. Indeed, if $\phi(z)=\phi(w)$, we have 

\[
\left(\frac{z_1}{z_n},\dots,\frac{z_{n-1}}{z_n}\right)=\left(\frac{w_1}{w_n},\dots,\frac{w_{n-1}}{w_n}\right)
\]  

\vskip 2mm

We have $|z_n|=|w_n|=1$, so we put $\frac{z_n}{w_n}=e^{i\alpha}$ and we have $z=e^{i\alpha}w$.

\vskip 2mm

As a consequence, $\phi$ induces a map $\overline{\phi}:\cN\rightarrow \cZ_{\cP,n-1}$ which is continuous and bijective. Actually, this is an homeomorphism since the inverse map $\phi^{-1}$ is the continuous map

\[
\begin{array}{rcccl}
\phi^{-1} & : & \cZ_{\cP,n-1} & \rightarrow & \cN \\
	       	&   & z	 				 		& \mapsto 	  & [(z,1)] 
\end{array}
\]

where $[(z,1)]$ denotes the equivalence class of $(z,1)\in \cZ_{\cP,n}$ of the diagonal action. 
\end{demo}

\vskip 2mm

\begin{cor}\label{complexe3} Let $\cP$ be a rationally starshaped sphere. Then there exists $N\in\N^*$  such that $\cZ_{\cP,N}$ can be endowed with a complex structure as LVMB manifold. 
\end{cor}

\begin{demo} Since $\cP$ is a rationally starshaped sphere, there exists a good system $(\cE,l)$ with type $(2m+1,n,k)$ whose associated complex is $\cP$ (cf. the previous subsection). Moreover, our construction of $(\cE,l)$ implies that $k>0$. So, by proposition \ref{complexe2}, $\cN$ is homeomorphic to $\cZ_{\cP,n-1}$. So, we put $N=n-1$ and we endow $\cZ_{\cP,N}$ with the complex structure induced by this homeomorphism. 
\end{demo}

\vskip 2mm

\begin{rmq} Let $N_0$ be the smallest integer $N$ as in corollary \ref{complexe3}. Then, for every $q\in\N$, $\cZ_{\cP,N_0+2q}$ can also be endowed with a complex structure and then we have $\cZ_{\cP,N_0+2q}=\cZ_{\cP,N_0}\times(S^1)^{2q}$. Indeed, let $\Lambda$ be the matrix whose columns are the vectors of the good system $(\cE,l)$ constructed in the proof of corollary \ref{complexe3}. We put $\widetilde{\cE}=\{\ P\cup\{N_0+1,N_0+2\}/\ P \in \cE \}$ and we define $\lambda_1,\dots,\lambda_{n+2}$ as the columns of the matrix 

\[
\begin{pmatrix}
\Lambda 		& 0  & 0 \\
-1 \dots -1 & 1  & 0 \\
-1 \dots -1 & -1 & 1 
\end{pmatrix}
\]

\vskip 2mm

Then it is easy to show that $(\widetilde{\cE},(\lambda_1,\dots,\lambda_{n+2}))$ is a good system and 

\[
\cZ_{\cP,N_0+2}=\cZ_{\cP,N_0}\times(S^1)^2
\]

\end{rmq}

\bibliographystyle{short}
\bibliography{versionanglaise}

\end{document}